\documentstyle[11pt]{article}

\begin{document}

\title{Equations of low-degree Projective Surfaces
with three-divisible Sets of Cusps}
\author{Wolf P. Barth \and S{\l}awomir Rams}
\date{{}}
\maketitle

\begin{abstract}
We determine the equations of surfaces
$Y \subset \mbox{{\bf P}}_3(\mbox{{\bf C}})$
of degrees
$\leq 6$ carrying a minimal, non-empty, three-divisible set of 
cusps.
\end{abstract}

\def\thefootnote{}
\footnotetext{
\hspace{-5.5ex} 
Suported by the DFG Schwerpunktprogramm
"Global methods in complex geometry". The second
author is supported by a
Fellowship of the Foundation for Polish Science
and KBN Grant No. 2 P03A 016 25. }

%


\def\Rn{\R^n}
\def\Rm{\R^m}
\def\Rz{\R^2}
\def\Rd{\R^3}

\def \qed {\hspace*{\fill}\frame{\rule[0pt]{0pt}{8pt}\rule[0pt]{8pt}
{0pt}}\par}
\def \qedup{\vskip-20pt\qed}

\def\R{{\rm I\!R}} 
\def\N{{\rm I\!N}} 
\def\F{{\rm I\!F}}
\def\M{{\rm I\!M}}
\def\H{{\rm I\!H}}
\def\K{{\rm I\!K}}
\def\P{{\rm I\!P}}
\def\E{{\mathchoice {\rm 1\mskip-4mu l} {\rm 1\mskip-4mu l}
{\rm 1\mskip-4.5mu l} {\rm 1\mskip-5mu l}}}
\def\C{{\mathchoice {\setbox0=\hbox{$\displaystyle\rm C$}\hbox{\hbox
to0pt{\kern0.4\wd0\vrule height0.9\ht0\hss}\box0}}
{\setbox0=\hbox{$\textstyle\rm C$}\hbox{\hbox
to0pt{\kern0.4\wd0\vrule height0.9\ht0\hss}\box0}}
{\setbox0=\hbox{$\scriptstyle\rm C$}\hbox{\hbox
to0pt{\kern0.4\wd0\vrule height0.9\ht0\hss}\box0}}
{\setbox0=\hbox{$\scriptscriptstyle\rm C$}\hbox{\hbox
to0pt{\kern0.4\wd0\vrule height0.9\ht0\hss}\box0}}}}
\def\bbbe{{\mathchoice {\setbox0=\hbox{\smalletextfont e}\hbox
{\raise
0.1\ht0\hbox to0pt{\kern0.4\wd0\vrule width0.3pt height0.7\ht0\hss}
\box0}}
{\setbox0=\hbox{\smalletextfont e}\hbox{\raise
0.1\ht0\hbox to0pt{\kern0.4\wd0\vrule width0.3pt height0.7\ht0\hss}
\box0}}
{\setbox0=\hbox{\smallescriptfont e}\hbox{\raise
0.1\ht0\hbox to0pt{\kern0.5\wd0\vrule width0.2pt height0.7\ht0\hss}
\box0}}
{\setbox0=\hbox{\smallescriptscriptfont e}\hbox{\raise
0.1\ht0\hbox to0pt{\kern0.4\wd0\vrule width0.2pt height0.7\ht0\hss}
\box0}}}}
\def\Q{{\mathchoice {\setbox0=\hbox{$\displaystyle\rm Q$}\hbox
{\raise
0.15\ht0\hbox to0pt{\kern0.4\wd0\vrule height0.8\ht0\hss}\box0}}
{\setbox0=\hbox{$\textstyle\rm Q$}\hbox{\raise
0.15\ht0\hbox to0pt{\kern0.4\wd0\vrule height0.8\ht0\hss}\box0}}
{\setbox0=\hbox{$\scriptstyle\rm Q$}\hbox{\raise
0.15\ht0\hbox to0pt{\kern0.4\wd0\vrule height0.7\ht0\hss}\box0}}
{\setbox0=\hbox{$\scriptscriptstyle\rm Q$}\hbox{\raise
0.15\ht0\hbox to0pt{\kern0.4\wd0\vrule height0.7\ht0\hss}\box0}}}}
\def\T{{\mathchoice {\setbox0=\hbox{$\displaystyle\rm
T$}\hbox{\hbox to0pt{\kern0.3\wd0\vrule height0.9\ht0\hss}\box0}}
{\setbox0=\hbox{$\textstyle\rm T$}\hbox{\hbox
to0pt{\kern0.3\wd0\vrule height0.9\ht0\hss}\box0}}
{\setbox0=\hbox{$\scriptstyle\rm T$}\hbox{\hbox
to0pt{\kern0.3\wd0\vrule height0.9\ht0\hss}\box0}}
{\setbox0=\hbox{$\scriptscriptstyle\rm T$}\hbox{\hbox
to0pt{\kern0.3\wd0\vrule height0.9\ht0\hss}\box0}}}}
\def\bbqS{{\mathchoice
{\setbox0=\hbox{$\displaystyle \rm S$}\hbox{\raise0.5\ht0\hbox
to0pt{\kern0.35\wd0\vrule height0.45\ht0\hss}\hbox
to0pt{\kern0.55\wd0\vrule height0.5\ht0\hss}\box0}}
{\setbox0=\hbox{$\textstyle \rm S$}\hbox{\raise0.5\ht0\hbox
to0pt{\kern0.35\wd0\vrule height0.45\ht0\hss}\hbox
to0pt{\kern0.55\wd0\vrule height0.5\ht0\hss}\box0}}
{\setbox0=\hbox{$\scriptstyle \rm S$}\hbox{\raise0.5\ht0\hbox
to0pt{\kern0.35\wd0\vrule height0.45\ht0\hss}\raise0.05\ht0\hbox
to0pt{\kern0.5\wd0\vrule height0.45\ht0\hss}\box0}}
{\setbox0=\hbox{$\scriptscriptstyle\rm S$}\hbox{\raise0.5\ht0\hbox
to0pt{\kern0.4\wd0\vrule height0.45\ht0\hss}\raise0.05\ht0\hbox
to0pt{\kern0.55\wd0\vrule height0.45\ht0\hss}\box0}}}}
\def\Z{{\mathchoice {\hbox{$\sf\textstyle Z\kern-0.4em Z$}}
{\hbox{$\sf\textstyle Z\kern-0.4em Z$}}
{\hbox{$\sf\scriptstyle Z\kern-0.3em Z$}}
{\hbox{$\sf\scriptscriptstyle Z\kern-0.2em Z$}}}}

\def\vecx{{\bf x}}
\def\vecxi{{\bf \xi}}
\def\vecy{{\bf y}}
\def\vecz{{\bf z}}
\def\norx{\parallel \vecx \parallel}
\def\nory{\parallel \vecy \parallel}
\def\norh{\parallel \vech \parallel}
\def\vecv{{\bf v}}
\def\vecw{{\bf w}}
\def\norv{\parallel \vecv \parallel}
\def\norw{\parallel \vecw \parallel}
\def\veca{{\bf a}}
\def\vecb{{\bf b}}
\def\nora{\parallel \veca \parallel}
\def\norb{\parallel \vecb \parallel}
\def\vecc{{\bf c}}
\def\vecd{{\bf d}}
\def\vece{{\bf e}}
\def\vecf{{\bf f}}
\def\vecg{{\bf g}}
\def\vech{{\bf h}}
\def\vect{{\bf t}}
\def\vecu{{\bf u}}
\def\vec0{{\bf 0}}
\def\veck{{\bf k}}
\def\vecm{{\bf m}}
\def\vecn{{\bf n}}
\def\vecp{{\bf p}}
\def\vecq{{\bf q}}
\def\vecr{{\bf r}}
\def\vecs{{\bf s}}

\def\vecA{{\bf A}}
\def\vecB{{\bf B}}
\def\vecC{{\bf C}}
\def\vecN{{\bf N}}

\def\On{{\bf O(n)}}
\def\SOn{{\bf SO(n)}}
\def\Tn{{\bf T(n)}}
\def\h{\hspace{0.1cm}}
\def\GLn{{\bf GL(n)}}
\def\GLnK{{\bf GL(n,\K)}}
\def\GLnR{{\bf GL(n,\R)}}
\def\GLnC{{\bf GL(n,\C)}}
\def\Affn{{\bf Aff(n)}}
\def\AffnK{{\bf Aff(n,\K)}}

\def\grq{\grqq \hspace{1mm}}

\def\OO{{\cal O}}
\def\LL{{\cal L}}
\def\MM{{\cal M}}
\def\NN{{\cal N}}
\def\CC{{\cal C}}
\def\DD{{\cal D}}
\def\OF{{\cal F}}

\newcommand{\D}{\displaystyle}
\newcommand{\vv}{\vspace{0.3cm}}
\newcommand{\vV}{\vspace{0.5cm}}
\newcommand{\VV}{\vspace{2cm}}
\newcommand{\const}{\mbox{const}}
\def\ttt{\;^{\D t}\!\,}

\newtheorem{defi}{Definition}[section]
\newtheorem{prop}{Proposition}[section]
\newtheorem{theo}{Theorem}[section]
\newtheorem{lemm}{Lemma}[section]

\vV

\setcounter{section}{-1}

\section{Introduction}

We consider algebraic surfaces $Y \subset \P_3(\C)$. A {\em cusp}
(=singularity $A_2$) on $Y$ is a singularity near which
the surface is given in local (analytic) coordinates $x,y$ and $z$,
centered at the singularity,
by an equation
$$xy-z^3=0.$$
This is an isolated quotient singularity $\C^2/\Z_3$. A set 
$P_1,...,P_n$
of cusps on $Y$ is called $3$-divisible, if there is a cyclic global
triple cover of $Y$ branched precisely over these cusps.
Equivalently: If $\pi:X \to Y$ is the minimal desingularization
introducing two $(-2)$-curves $E_{\nu}', E_{\nu}''$ over each cusp,
then there is a way to label these curves such that the divisor 
class of
$$\sum_{\nu =1}^n (2E_{\nu}'+E_{\nu}'')$$
is divisible by 3 in $NS(X)$ \cite{t, b1}. {\em The aim of this 
note is to
determine the equations for surfaces $Y \subset \P_3$ of degrees 
$\leq 6$
carrying a minimal, non-empty set of 3-divisible cusps.}
Applications of 3-divisible sets
can be found in \cite{t}, \cite{kz}.

Recall \cite{t} that a non-empty 3-divisible set of cusps on a
surface of degree $d$ contains $n$ points with
$$ \begin{array}{c|ccc}
d & 3 & 4 & 5 \\ \hline
n & 3 & 6 & \geq 12 \\
\end{array} $$

Here we show

\vV
\noindent
{\bf Theorem~\ref{atleasteighteen}}
{\em Any non-empty 3-divisible set of cusps on a sextic surface
contains at least 18 points.}

\vv
The cubic surface $Y \subset \P_3$ with three cusps is unique, up to
a choice of coordinates it has the equation $x_1x_2x_3-x_0^3=0$
\cite{t}.
For surfaces $Y \subset \P_3$ of degrees d=4, 5, 6 with
a 3-divisible set of
6, 12, 18 cusps we show (see lemma~\ref{quarticsandquintics}
and theorem~\ref{typesofsextics}):

\vV
\noindent
{\bf Theorem.} {\em There are two polynomials $s',s''$ of degree 3 
and some
polynomial $s$ of degree 2 such that the sextic polynomial
$s' \cdot s'' -s^3$ vanishes on $Y$. The equation $s' \cdot s''-
s^3=0$
in the case of degree

\begin{itemize}
\item[d=4:] can be chosen such that it
defines the quartic $Y$ together with some residual quadric;
\item[d=5:] defines the quintic $Y$ together with a residual plane;
\item[d=6:] either defines the sextic $Y$ or it vanishes 
identically.
If the latter holds then
$Y$ has an equation $l' \cdot l'' \cdot f -g^3=0$ with
$\deg(l')=\deg(l'')=1, \, \deg(g)=2, \, \deg(f)=4$.
\end{itemize} }

\vv
This description allows to conclude:

\begin{itemize}
\item Quartic surfaces with six irreducible cusps form an 
irreducible family.
Based on the equation given above, this has been shown in \cite{r}.
\item Quintic surfaces with 12 three-divisible cusps form
an irreducible family (thm.2.1).
\item Sextic surfaces with 18 three-divisible cusps form two 
irreducible
families. In fact the equations
for the surfaces in the two families are of
the two types given in the theorem above.
That these two families are disjoint,
this follows quite easily when considering the codes of these 
surfacesi \cite{br}. 
\end{itemize}

\vv
{\em Convention:} In this note the base field always is $\C$. When
we consider singular surfaces, by a divisor we mean a Weil-divisor.

\section{Contact cubics} \label{contactsurfaces}

Let $Y \subset \P_3$ be a surface of degree $d=4,5$ or $6$. Let
$P_1,...,P_n \in Y$ be a 3-divisible set of $n=6,12$ or $\leq 18$
cusps. We assume that $Y$ is smooth but for these cusps,
and perhaps for some further rational double points. 
The aim of this section is to show that there are two cubic surfaces 
touching $Y$ to the third order along two curves passing
through the n cusps.

We denote by
$$\pi: X \to Y$$
the minimal resolution of singularities. In this situation it is
well-known, that the canonical bundle of $X$ is
$$K_X = \pi^*( \OO_Y(d-4)).$$
By abuse of notation we write $\OO_X(k):=\pi^*\OO_Y(k), \, k \in 
\Z$.

The minimal desingularisation introduces two $(-2)$-curves
$E_{\nu}'$ and $E_{\nu}''$ over each cusp $P_{\nu}$. We label them 
such that
the two classes
$$\LL':=\frac{1}{3} \sum_{\nu=1}^n (2E_{\nu}'+E_{\nu}'')
\quad \mbox{ and } \quad
\LL'':=\frac{1}{3}\sum_{\nu=1}^n (E_{\nu}'+2 E_{\nu}'')$$
exist in $NS(X)$. We also introduce the divisor classes
$$ \CC':=\OO_X(1)-\LL', \, \, \CC'':=\OO_X(1)-\LL'', \quad
\DD':=\OO_X(1)+\LL', \, \, \DD'':=\OO_X(1)+\LL'' $$
on $X$.
Obviously
$$(\LL' . E_{\nu}')= -1, \quad (\LL'. E_{\nu}'')=0$$
and similarly for $\LL''$. This implies
$$(\LL')^2 = (\LL'')^2= - \frac{2}{3}n, \quad
(\CC')^2=(\CC'')^2=(\DD')^2=(\DD'')^2= d-\frac{2}{3}n.$$
The assertions stated for $\LL'$ in the next lemma
similarly hold for $\LL''$.

\begin{lemm} \label{lemmaoneone}
{\bf a)} The class $\LL'$ in $NS(X)$
is not an integral linear combination
$$E+\sum_{\nu=1}^n (a_{\nu}'E_{\nu}' + a_{\nu}'' E_{\nu}''), \quad
a_{\nu}', a_{\nu}'' \in \Z,$$
with $E$ an exceptional divisor lying over the additional 
singularities
of $Y$.

{\bf b)} The class $\LL'$ is not effective.

{\bf c)} [R 2, lemma~2.1] If $C' \in |\CC'|$, then $\pi(C') \subset 
Y$
cannot be a plane section.
\end{lemm}

Proof. a) If $\LL'$ is represented by an integral linear combination
$L$ as in the assertion, we find
$$-1=(L.E_{\nu}')=a_{\nu}''-2a_{\nu}', \quad
0=(L.E_{\nu}'') = a_{\nu}'-2a_{\nu}''.$$
This would imply $a_{\nu}'=2a_{\nu}''$ and lead to
the contradiction $a_{\nu}''=\frac{1}{3}$.

b) From $(\OO_X(1).E_{\nu}')=(\OO_X(1).E_{\nu}'')=0$ we conclude
$(\OO_X(1).\LL')=0$. If $L$ is an effective divisor representing 
$\LL'$,
then $(\OO_X(1).L)=0$ implies that $L$ is an integral linear 
combination
of exceptional divisors lying over the singularities of $Y$. This is
in conflict with a).

c) Assume that $\pi(C')$ is a plane section. The total transform
of this plane section on $X$ then is a divisor
$C'+E$ with $E$ an integral linear combination of exceptional 
divisors.
Hence
$$\OO_X(1) \sim C'+E \sim \OO_X(1)-\LL'+E,
\quad \mbox{ so } \quad \LL' \sim E,$$
a contradiction with a).
\qed

\vv
We use Riemann-Roch on $X$
$$\chi(\CC')=\chi(\OO_X) + \frac{1}{2}\CC' . (\CC' - K_X)$$
and similarly for $\CC''$. Here $\chi(\OO_X)$ is the
Euler-characteristic $\chi(\OO_{Y'})$ for any smooth
surface $Y'$ of degree $d$, since $X$ and $Y'$ are
diffeomorphic \cite{brie}.
We obtain the table
$$ \begin{array}{l|cccc|c}
d & K_X & \chi(\OO_X) & (\CC')^2 & -\CC'. K_X & \chi(\CC') \\ \hline
4 & \OO_X & 2 & 0 & 0 & 2\\
5 & \OO_X(1) & 5 &-3 & -5 & 1 \\
6 & \OO_X(2) & 11 & 6 - \frac{2}{3}n & -12 & 8 - \frac{n}{3} \\
\end{array} $$

The main fact we need is this:

\begin{prop} \label{classesareeffective}
The divisor classes $\CC'$ and $\CC''$ on $X$ are effective.
\end{prop}

Proof. It suffices to prove the assertion for $\CC'$. This follows 
from
Riemann-Roch by controlling $h^2(\CC')$.
The proof is easy for $d=4$ or $5$, but quite tedious for $d=6$.
We consider these three cases.

$d=4$: Here $h^2(\CC')=h^0(-\CC')=h^0(\OO_X(-1)+\LL')$. Since $\OO_X
(1)$
is nef, from $(-\CC'.\OO_X(1))=-4$ it follows that $-\CC'$ cannot 
be effective.
So $h^2(\CC')=0$ and $h^0(\CC') \geq 2$.

$d=5$: Now $h^2(\CC')=h^0(\OO_X(1)-\CC')=h^0(\LL')=0$
by lemma~\ref{lemmaoneone}~b).

\vv
So let us now concentrate on the case $d=6$. Since $K_X=\OO_X(2)$,
by Serre-duality
$$h^2(\CC')=h^0(\DD'), \quad h^2(\CC'')=h^0(\DD'').$$
The effective class $E:=\sum (E_{\nu}'+E_{\nu}'')$ defines
a map $\CC' \to \DD''$ and an exact sequenc
$$0 \to \CC' \to \DD'' \to \DD''|E \to 0.$$
From $(\DD''.E_{\nu}')=0, \, (\DD''.E_{\nu}'')=-1$ we
conclude $h^0(\DD''|E)=0$. This implies
$$h^0(\CC')=h^0(\DD'')=h^2(\CC'').$$
So proposition~\ref{classesareeffective} follows from Riemann-Roch 
even in the case
$d=6$ and $n \leq 18$ if we show

\begin{prop} \label{nomorethanone}
If $d=6$, then $h^0(\CC') \leq 1, \, h^0(\CC'') \leq 1$.
\end{prop}

Proof. As usual this statement needs to be proven for $\CC'$ only. 
So let us
assume the assertion to be false, and $h^0(\CC') \geq 2$. We 
consider
the linear system $|\CC'|$ on $X$ with
$$\CC'.\OO_X(1) = 6, \quad (\CC')^2=6+(\LL')^2 \geq -6.$$
Let $B$ be the base curve of the system and $|F|$ its free part.
Then $F^2 \geq 0$ and
$$ 1 \leq (F.\OO_X(1)) \leq (\CC.\OO_X(1)) =6.$$

\begin{lemm} We have $F^2=0,2$ or $4$. If $F^2=2$, then
$(F.\OO_X(1)) \geq 4$ and if $F^2=4$, then $(F.\OO_X(1)) \geq 5$.
\end{lemm}

Proof. Since $K_X=\OO_X(2)$ is a square, from Riemann Roch it 
follows that
all divisor classes on $X$ have an even self-intersection.
We consider the determinant
$$ \det \left( \begin{array}{cc}
\OO_X(1)^2 & (F.\OO_X(1)) \\
(F.\OO_X(1)) & F^2 \\ \end{array} \right) =
6 \cdot F^2-(F.\OO_X(1))^2.$$
Since $(F.\OO_X(1)) \leq 6$, the Hodge index theorem implies $F^2 
\leq 6$.
And if $F^2=6$, then $F \sim r \OO_X(1)$ with some rational number 
$r$.
From $6r = (F.\OO_X(1))=6$ we conclude $r=1$ and $F \sim \OO_X(1)$. 
But
this implies
$$\OO_X(1) - \LL' = \OO_X(1)+B$$
and $-\LL' \sim B$ would be effective. Since $3\LL'$ is effective, 
this is
a contradiction.

We are left with the cases $F^2=0,2$ or $4$. In these cases the 
class of $F$
cannot be a rational multiple of $\OO_X(1)$. By the Hodge index
theorem the above determinant must be negative. This implies
$(F.\OO_X(1)) \geq 4$ for $F^2=2$ and $\geq 5$ for $F^2=4$. \qed

\begin{lemm} \label{thecurveisirreducible}
The general curve $F_0 \in |F|$ is irreducible.
\end{lemm}

Proof. Assume that all curves $F_0 \in |F|$ are reducible. First
we show that then $F \sim k A, \, 2 \leq k \in \N,$ with some 
effective
irreducible curve $A \subset X$. To do this, blow up the base 
points of $|F|$
via $\sigma:X' \to X$. (If $F^2=0$ there are no base points.)
Let $|F'|$ on $X'$ be the proper transform of the linear system 
$|F|$.
If $(F')^2>0$, by Bertini the general curve $F_0' \in |F'|$ will be
irreducible, and this implies that the general curve $F_0$ is
irreducible. So we are in the case $(F')^2=0$. Then $F'$ is composed
of a pencil. All irreducible components $A_i$ of a general fibre 
$F_0'$
are algebraically equivalent, and by $q(X')=0$ linearly equivalent.
This implies $F_0 \sim k A'$ with $A'$ such an irreducible component
and $k$ the number of components in a fibre. But with $A:=\sigma(A')
$
we get
$$F_0 \sim \sigma(F_0') \sim k A.$$

Now consider the different cases: If $F^2=2$, then $F$ cannot
be linearly equivalent to some $k \cdot A$, because then $F^2\geq 4 
A^2$.
And if $F^2=4$, then $F^2 = k^2 A^2$ would imply $k=2$ and
$A^2=1$. This is impossible, because all self-intersections of 
curves on $X$
are even numbers. We are left with the case $F^2=0$. Here
$F \sim k \cdot A$ with a {\em smooth irreducible}
curve $A \subset X$. The degree
of its (birational) image $\pi(A) \subset \P_3$ can be at most
$6/2=3$. So $A$
will be rational or elliptic. But this is in conflict with the
adjunction formula
$$2g(A)-2 = 2 (A.\OO_X(1))$$
and $(A.\OO_X(1)) \geq 1$. \qed

\begin{lemm} The image curve $\pi(A)$ of a general $A \in |F|$ is a 
plane
quintic. \end{lemm}

Proof.
Put $a:=(A.\OO_X(1))$. Here $a > 0$, since $A$ is not exceptional.
It suffices to prove that $\pi(A)$ is planar irreducible.
The adjunction formula for the plane curve $\pi(A)$ then shows
$$A^2+2a =2p_a(A)-2 \leq a(a-3).$$
From $A^2 \geq 0$ we conclude $a^2-5a \geq 0$. So $a \geq 5$
and $a=6$ only if $A$ is a plane section of $Y$. But the latter
would contradict lemma~\ref{lemmaoneone}~c).

If $F^2=0$,
the general fibre $A$ of $F$ is smooth by Bertini, since the system
$|F|$ has no base point. If $F^2=2$ or $4$, by Bertini the general
fibre $A$ can be singular only in base points. If all fibres $A$
are singular in some base point, then necessarily $A^2 \geq 4$.
In this case the base point is a double point on all fibres. Two 
distinct
fibres cannot have a common tangent at this singularity. This 
implies
that the general curve $A$ has an ordinary node at the base point 
and is smooth
in all other points. The curve $\pi(A)$ is irreducible of degree 
$a$.

We apply the adjunction formula on $X$:
$$2 p_a(A)-2 = A^2+2a, \quad p_a(A) =a+\frac{1}{2}A^2+1.$$
If $A^2=0$ or $2$, the general curve $A$ is smooth with $p_g(A)=p_a
(A)$.
For $A^2=4$ the general curve $A$ may have an ordinary node, hence
$p_g(A) \geq p_a(A)-1$. In these three cases we have the following
lower bounds for $p_g(A)$:
$$ \begin{array}{c|ccc}
A^2 & 0 & 2 & 4 \\ \hline
p_g(A) & a+1 & a+2 & a+2 \\ \end{array} $$
Since $a \geq 1$, necessarily $p_g(A) \geq 2$. All irreducible
space curves of degrees $\leq 4$ are rational or elliptic, unless
they are plane quartics, in which case we are done. This already 
shows
$a \geq 5$. We apply the Castelnuovo-bound
[ACGH, p.~116] for the irreducible curve $\pi(A)$,
if it is non-planar:
$$ \begin{array}{c|ccc}
a & 5 & 6 \\ \hline
p_g(A) \leq & 2 & 4 \\
\end{array}$$
It cannot be met by our curve, and this curve must be planar.
\qed

\vv
Now we aim at the final contradiction. Consider some curve
$B+F \in |\CC'|$ with $\pi(F)$ an irreducible planar quintic.
The plane $S_1$ of $\pi(F)$ cannot touch $Y$ along $\pi(F)$.
The curve $3 \cdot (B+F)$ has the class
$$\OO_X(3) -\sum_{\nu=1}^n (2E_{\nu}'+E_{\nu}'').$$
So there will be a cubic surface $S \subset \P_3$ pulling
back on $X$ to the curve $3(B+F)+\sum(2E_{\nu}'+E_{\nu}'')$.
In particular $S$ will meet the plane $S_1$ of $\pi(F)$ in this 
quintic
curve. Hence $S=S_1+S_2$ splits off $S_1$. The residual quadric
$S_2 \subset S$
will touch the surface $Y$ along the curve $\pi(F)$. So even $S_2$ 
splits
off the plane $S_1$, say $S_2=S_1+S_3$ with another plane $S_3$.
And again $S_3$ will contain the plane quintic, hence it coincides 
with $S_1$.
We have shown $S=3S_1$. Then $S_1$ must contain the whole curve
$\pi(B+F)$. This means that $\pi(B+F)$ is a plane section of $Y$,
contradicting lemma~\ref{lemmaoneone}~c).
This proves prop.~\ref{nomorethanone}. \qed

\vv
\begin{theo} \label{atleasteighteen}
Each non-empty 3-divisible set of cusps on a sextic surface $Y$
(with at most rational double points) contains at least $18$ points.
\end{theo}

Proof. By Riemann-Roch
$$2 = h^0(\CC')+h^0(\CC'') \geq \chi(\CC')=8 - \frac{n}{3}. \Box $$

\vv
\begin{theo} Assume that the surface $Y \subset \P_3$ of degree
$4,5$ or $6$ carries a 3-divisible set of cusps of the minimal
number $6,12$ or $18$. Then

\begin{itemize}
\item there are two distinct cubic
surfaces $S'$ and $S''$ passing through these cusps and touching $Y$
to the third order away from these cusps along two curves $\pi(C')$
and $\pi(C'')$;
\item there is a quadric surface $S$ cutting out on $Y$
the two curves
$\pi(C')$ and $\pi(C'')$. \end{itemize} \end{theo}

Proof. By prop.~\ref{classesareeffective} the classes $\CC'$ and 
$\CC''$ are effective.
Choose divisors $C' \in |\CC'|, \, C'' \in |\CC''|$. Then
$$3C' \sim \OO_X(3)- \sum (2E_{\nu}'+E_{\nu}''), \quad
3C'' \sim \OO_X(3)- \sum (E_{\nu}'+2E_{\nu}'').$$
This implies that there are cubic surfaces $S'$ and $S''$
cutting out on $Y$ the divisors $3 \pi(C')$ and $3 \pi(C'')$
respectively. If these surfaces would coincide, then so
would their divisor $\pi(C')=\pi(C'')$ of contact with $Y$.
Let $C \subset X$ be the proper transform of this divisor. Then
$$C'=C+E', \quad C'' = C+E''$$
with integral linear combinations $E',E''$ of exceptional divisors.
But then
$$C'-C'' \sim \frac{1}{3} \sum (-2E_{\nu}'-E_{\nu}''+E_{\nu}'+2E_
{\nu}'')
=\frac{1}{3} \sum (E_{\nu}''-E_{\nu}') \sim \LL'-\sum E_{\nu}'$$
would be linearly equivalent to $E'-E''$,
a contradiction with lemma~\ref{lemmaoneone}~a).

The divisor $C'+C''$ on $X$ has the class
$$\OO_X(2)-\sum (E_{\nu}'+E_{\nu}'').$$
Hence $\pi(C'+C'')$ is the intersection of $Y$
with a quadric surface $S$.
\qed

\vV

\section{Equations}

Consider the surface $Y$ of degree 4, 5 or 6 with its contact
cubics
$$S':\, s'=0, \quad S'': \, s''=0$$
as constructed at the end of section~\ref{contactsurfaces}.
The product $s' \cdot s''$
vanishes on the divisor $\pi(C'+C'')$ cut out on $Y$ by some
quadric surface
$$S: \, s=0.$$
In particular $s' \cdot s''/s^3$
is constant on $Y$. This constant may be squeezed e.g. into $s^3$. 
Then the
sextic polynomial
$$s' \cdot s'' - s^3$$
vanishes on $Y$. Unless it vanishes identically on $\P_3$, this 
gives
a sextic equation $s' \cdot s'' -s^3=0$ which defines in degree

\begin{center} \begin{tabular}{ll}
4 & the surface $Y$ together with a residual quadric, \\
5 & the surface $Y$ together with a residual plane, \\
6 & the surface $Y$ properly. \end{tabular} \end{center}

We now analyze the situation
$$s' \cdot s'' -s^3 = 0 \quad \mbox{identically on } \P_3.$$
If here $s$ would be irreducible, it would divide $s'$ and $s''$, 
say
$$s'=s \cdot l', \quad s''=s\cdot l'' \quad \mbox{with} \quad
\deg(l')=\deg(l'')=1.$$
But then $l' \cdot l''$ would divide $s$, a contradiction.
So necessarily $s=l' \cdot l''$ is reducible and
$s'\cdot s'' = (l' \cdot l'')^3$. Here $s' = \const \cdot (l')^3$
or $s'' = \const \cdot (l'')^3$ would contradict lemma~\ref
{lemmaoneone}~b).
W.l.o.g.
we therefore may assume
$$s'=(l')^2 \cdot l'', \quad s''=l' \cdot (l'')^2, \quad
l',l'' \mbox{ linearly independent}.$$

Let $H'$ and $H'' \subset \P_3$ be the planes of equation
$l'=0, \, l''=0$ respectively. Let $I=H' \cap H''$ be their
line of intersection. Consider the divisors cut out on $Y$
by these two planes. If $Y$ contains the line $I$, let
$i'$, resp. $i''$ be the multiplicities of $I$ in these
divisors. So
$$H' . Y = A' + i'I, \quad H'' . Y = A'' + i'' I,$$
with $I \not\subset A',A''$ and
$i'=i''=0$ if $I \not\subset Y$.
Recall from section~\ref{contactsurfaces}
$$ 2A'+(2i'+i'')I+A''=S'. Y = 3 \pi(C'), \quad
A'+(i'+2i'')I+2A''=S''. Y = 3 \pi(C'').$$
This implies $A'=3B'$ and $A''=3B''$. In particular
both the divisors $\pi(C')$ and $\pi(C'')$ will be non-reduced.

\begin{lemm} \label{quarticsandquintics}
{\bf a)} If $d=4$, the divisors $C'$ and $C''$ can be chosen
such that the sextic $s' \cdot s'' - s^3$ will not vanish 
identically.

{\bf b)} If $d=5$ the situation that $s' \cdot s''-s^3$ vanishes
identically does not occur. \end{lemm}

Proof. a) If $d=4$ we know
$$\chi(\CC')=\chi(\CC'')= 2.$$
For $C' \in |\CC'|$ the image $\pi(C')$ has degree 4. If $C'$ is
chosen generically, this image cannot consist of rational curves
only. It is an elliptic quartic curve, or an elliptic cubic plane
curve together with a line. In particular it is reduced. Hence
$s' \cdot s'' - s^3$ will not vanish identically on $\P_3$.

b) If $d=5$, both the divisors
$$Y . H' = 3B'+i'I, \quad Y . H''=3B''+i''I$$
have degree 5. This implies $i' \geq 2$ and $i'' \geq 2$. But
this means, the surface $Y$ touches both the planes $H'$
and $H''$ along $I$. Then $I$ will be a double line on $Y$.
Since $Y$ has isolated singularities only, this is impossible.
\qed

\vv
\begin{theo} Quintics with 12 three-divisible cusps (and
else at most rational double points) form an irreducible family.
\end{theo}

Proof. Let $Y$ be a quintic with 12 three-divisible cusps.
Choose polynomials $s',s'',s$ as in thm. 1.2 such that 
$s' \cdot s'' -s^3$ vanishes on $Y$ together with some residual
plane $R$. By equivariance under the projective group it suffices to
show that the equations $s' \cdot s'' -s^3$ vanishing on a fixed
plane $R$ form an irreducible family. So fix the plane
$R: \, x_0=0$. We have to consider separately different (mostly
degenerate) cases:

{\em Case 1:} The quadric $S: \, s=0$ contains the plane $R$.

{\em Subcase 1a)} $R=S$ and, say, $s=x_0^2, \, s'=x_0 \cdot q'$.
Then $Y$ has an equation $f=q' \cdot s'' - x_0^5=0$ with all 
the twelve cusps on $R$. The set $q'=s''=x_0=0$ consists
of singularities of $Y$, hence this set is discrete. In particular
it contains at most six points, which can be cusps of $Y$.
So there must be at least six more cusps on $Y$ where $s''=0$
but $q' \not= 0$. But there $d(q' \cdot s'')=q' \cdot ds''$
vanishes if and only if $ds''=0$, i.e., if the cubic surface
$S''$ itself is singular. If the curve $s''=x_0=0$ is reduced,
this can happen in at most three points. So the curve
$s''=x_0=0$ must be non-reduced. It cannot consist of a threefold
line, because on this line there could be at most four singularities
of $Y$. And if it consists of a double line together with
a simple line, on the simple line there cannot be a singularity of $S''$.
Again there are at most four singularities of $Y$ on the double line,
not enough. This excludes subcase 1a). 

{\em Subcase 1b)} The quadric $S$ splits off $R$, say
$s=x_0 \cdot x_1$ and $f=q' \cdot s'' -x_0^2 \cdot x_1^3$. Denote by
$P$ the plane $x_1=0$ and by $L$ the line $P \cap R$. Again the set
$q'=s''=x_1=0$ is singular on $Y$, hence discrete, and it
can contain at most six cusps of $Y$. At least six more cusps of $Y$
must lie on the curve $C: \, s''=x_0=0$ away from $q'|L=0$.  Along this
curve $C$, away from $L$ and $q'=0$ locally $x_0^2|Y \sim s''|Y$
is a third power. This implies $x_0|Y$ itself is a third power. Hence
in those points $S''$ and $R$ touch. For the intersection
curve $S'' . R$ the possibilities are ($M \subset R$ a line $\not=L$):

\begin{itemize}
\item $S'' . R = 3M$: On $M$ there are at most four singularities of
$Y$ with no further singularities on $L$.
\item $S''=2M+L$: Again on $M$ there are at most four singularities of $Y$. Away
from $M$ on $L$ we have $ds'' \not= 0$ and singularities could occur
only in points where $q'=0$. But those points belong to the 
$\leq 6$ points in $q'=s''=x_1=0$ considered already.
\item $S'' . R = 3L$: Again there are at most four singularities of
$Y$ on this line.
\end{itemize}
In all these cases there are less than 12 cusps. 

\vv
{\em Case 2:} The quadric $S$ meets $R$ in two distinct lines, say
$L_1:x_0=x_1=0$ and $L_2:x_0=x_2=0$. So $s=x_1 \cdot x_2+x_0 \cdot l$.

{\em Subcase 2a)} $s'|R=x_1^3|R$ and $s''|R=x_2^3|R$. This is the general
case. It is analyzed in [BR, thm.1.1] There it is shown that for general choice
of the equation the surface $Y$ indeed has twelve cusps away from $R$
and is smooth in all its other points. That surfaces with such an
equation form an irreducible family, this is quite obvious.

{\em Subcase 2b)}: $s'|R=x_1^2 \cdot x_2|R$ and $s''|R=x_1 \cdot x_2^2|R$.
Then $s'=x_1^2 \cdot x_2+x_0 \cdot q', \, 
s'' = x_1 \cdot x_2^2  +x_0 \cdot q''$ and $Y$ is defined by
$$f:=x_1^2 x_2 q''+ x_1 x_2^2 q' + x_0 q' q''
-3x_1^2x_2^2l-3x_0x_1x_2l^2-x_0^2l^3.$$
This equation vanishes on both the lines $L_i$, hence both lines
belong to $Y$. We claim $q'|L_1=0$ (and similarly 
$q''|L_2=0$). Indeed, this is true if $S'$ is singular along $L_1$.
And if $S'$ is smooth along $L_1$ it touches there $Y$ and $R$
simultaneously. This implies that $Y$ and $R$ touch along $L_1$.
And by the form of $f$, also this proves the claim. So in fact
$Y$ is singular along both lines, contradiction.

\vv
{\em Case 3:} The quadric $S$ meets $R$ in a repeated line, say
$s=x_1^2+x_0 \cdot l$ and and $s'=x_1^3+x_0 \cdot q', \,
s''= x_1^3 +x_0 \cdot q''$. Then $Y$ is defined by the polynomial
$$f:=x_1^3(q'+q'')+x_0 q' q'' -3x_1^4 l - 3x_0x_1^2 l^2 -x_0^2 l^3.$$
By letting the lines $L_1$ and $L_2$ in subcase 2a) move together 
one proves that the surfaces in this family lie in the closure of
the irreducible family from subcase 2a). However there is still some subtile
point, namely, whether the twelve cusps from 2a) converge to the twelve
3-divisible cusps on the surfaces from case 3. Therefore let us
analyse the possible configuration of cusps in case 3:

By the form of $f$ the line $L: \, x_0=x_1=0$ lies on $Y$. Along this line
$S'$ and $S$, as well as $S''$ and $S$ touch. So off $R$ the (there necessarily
discrete) set $S \cap S' \cap S''$ contains at most eight points.
And they are cusps on $Y$ if and only if this intersection is transversal,
i.e., if it are eight points indeed. Additional cusps of $Y$ can lie
only on $S' \cap S'' \cap R = L$. Since $(\partial_0 f)|L=q' \cdot q''$
there are indeed four more singularities of $Y$ on this line,
for general choice of $q'$ and $q''$. That these are cusps, 
e.g. the points of $L$ where $q'=0$, this
follows from the SQH-criterion [BW] with the weights
$$wt(x_1)=\frac{1}{3}, \quad wt(x_0)=wt(q')=\frac{1}{2}.$$

It remains to show for a family $Y_t$ of surfaces from
subcase 2a) converging to a surface $Y_0$ considered here, with
twelve cusps, that the twelve cusps on $Y_t$ converge point
by point to the twelve cusps on $Y_0$. If this would not be the case,
at least two cusps $P_1(t)$ and $P_2(t)$ on $Y_t$ would converge
to the same cusp $P_0$ on $Y_0$. Now fix a small sphere $S$
around $P_0$ containing $P_1(t)$ and $P_2(t)$ but no other
singularities of $Y_0$. We consider the Milnor number
$\mu(P_0)=2$ of $Y_0$ at $P_0$. By Milnor's definition
[Mi, p. 59] it is the degree of the map of $S$ onto the unit
sphere, defined by $g/\parallel g \parallel$ with
$g = (\partial_0 f_0, \partial_1 f_0, \partial_2 f_0)$, where
$f_0=0$ is a local equation for $Y_0$ and the derivatives
are taken w.r.t. some local coordinates. By continuity this is the same
as the number defined by $f_t$, where $f_t=0$ for small $t$ is a local
equation of $Y_t$, such that $f_t$ converges to $f_0$. But by
[Mi, p.112] this is at least the sum of the Milnor numbers of all singularities
of $Y_t$ contained in the fixed sphere. This sum is 
$\geq \mu(P_1(t))+\mu(P_2(t))=4$, contradiction. \qed 

\vv
Now we turn to degree 6. If $s' \cdot s''-s^3$ does not vanish
identically on $\P_3$, then $s' \cdot s''-s^3=0$ is an equation for
the sextic surface $Y$. However the other case is possible too:

\begin{theo} \label{typesofsextics}
Assume that the polynomial $s' \cdot s''-s^3$ vanishes
identically on $\P_3$. Then there are homogeneous polynomials
$l',l'',g,f$ of degrees
$$\deg(l')=\deg(l'')=1, \quad \deg(g)=2, \quad \deg(f)=4,$$
such that $Y$ is defined by the equation
$$l' \cdot l'' \cdot f - g^3=0.$$
\end{theo}

Proof. Put $s'=(l')^2 \cdot l''$ and $s''=l' \cdot (l'')^2$ as
above. We know already, that these polynomials define divisors
$$3B'+i'I \mbox{ on } H', \quad 3B''+i''I \mbox{ on } H''.$$
Here $i'=i''=0 \, \mbox{mod} \, 3$. If $i' \not= 0$, then the line
$I$ is contained in $Y$ and $Y$ touches $H'$ along this line.
But then $i'' \not= 0$ too, and $Y$ also touches $H'$ along $I$.
So $Y$ would contain $I$ as a multiple curve, which was
excluded. In this way we see $i'=i''=0$.

Next we show that both the curves $B'$ and $B''$ are reduced
conics. Assume e.g. that $B'$ is a double line. All 18 cusps of
$Y$ lie on the union $H' \cup H''$ of both the planes, i.e. on
$B' \cup B''$. But since $Y$ is not singular along a curve, it
can have at most five cusps on the line $B'$ and at most ten
cusps on the conic $B''$, not enough.

Now we claim that there is a quadratic polynomial $g$ restricting to
$H'$ as an equation for $B'$ and to $H''$ as an equation for $B''$.
First observe that both the conics $B'$ and $B''$ meet $I$ in the
same point set $Y \cap I$. This set may consist of two distinct 
points or
of one point counted twice. Let $b'$ and $b''$ be two quadratic 
polynomials
restricting to the planes $H'$ and $H''$ as equations for $B'$ and 
$B''$
respectively. Then $b'|I$ and $b''|I$ differ by a constant factor.
We may adjust one of the quadric polynomials such that
$$b'|I=b''|I.$$
This implies
$$b'|H'' = b''|H''+(l' \cdot h)|H'', \quad \deg(h)=1.$$
Then the quadratic polynomial
$$g:=b'-l'h$$
restricts to $H'$ as $b'|H'$ and to $H''$ as $b''|H''$. In 
particular
the quadric $g=0$ cuts out on $H'$ the conic $B'$ and on $H''$
the conic $B''$.

Let now $\varphi=0$ be an equation for $Y$. Then both the 
polynomials
$\varphi$ and $g^3$ define on $H \cup H'$ the divisor $3(B'+B'')$.
On the line $I$ both polynomials differ by a constant factor.
We may adjust $\varphi$ by this factor, to obtain
$\varphi|I=g^3|I$. This then implies that $\varphi=g^3$ on
$H' \cup H''$. So there is a quartic polynomial $f$ with
$\varphi - g^3=l' \cdot l'' \cdot f$, or
$$\varphi = l' \cdot l'' \cdot f -(-g)^3$$
as we claimed. \qed

It is a different question, whether sextic surfaces with 18
three-divisible cusps and of equations
$$s' \cdot s'' - s^3=0, \quad \mbox{resp.} \quad
l' \cdot l'' \cdot f - g^3$$
with
$$\deg(l')=\deg(l'')=1, \quad \deg(s)=\deg(g)=2, \quad
\deg(s')=\deg(s'')=3, \quad \deg(f)=4$$
exist indeed. But choose all these polynomials generically.
Then a Bertini-type argument \cite{br} shows that the
surfaces with these equations are smooth but for the
$3 \cdot 3 \cdot 2 = 18$ intersection points
$$s'=s''=s=0$$
in the first case, and the $2+4 \cdot 2+4 \cdot 2=18$ points
$$ l'=l''=g=0, \quad l'=f=g=0, \quad l''=f=g=0$$
in the second case. We refer to \cite{br} for the fact, that
these 18 points indeed form a 3-divisible set.

It is clear that the polynomials
$\varphi=s' \cdot s'' -s^3 \in H^0(\OO_{\P_3}(6))$ with
$s,s'' \in H^0(\OO_{\P_3}(3))$ and $s \in H^0(\OO_{\P_3}(2)$
varying arbitrarily form an irreducible variety. It contains of
course polynomials $\varphi$ defining degenerate surfaces,
namely

\begin{itemize}
\item $\varphi \equiv 0$,
\item defining surfaces with more than 18 singularities,
\item defining surfaces where the 18 singularities defined
degenerate or come together.
\end{itemize}

But these polynomials form a Zariski-closed subset. After removing 
them
one has an irreducible variety of sextic polynomials defining
surfaces $Y$ which are non-degenerate. This proves that
all non-degenerate sextic surfaces of equation $s' \cdot s'' -s^3=0$
as above form an irreducible family. The same argument of course
works for the second type of equation. We refer to \cite{br}
for the fact that these two families are disjoint.

\noindent
Wolf P. Barth. 
Mathematisches Institut der Universit\"{a}t, Bismarckstr. 1 1/2,
D - 91054 Erlangen, 
e-mail: barth@mi.uni-erlangen.de

\vspace{1ex}
\noindent
S{\l}awomir Rams.
Institute of Mathematics, Jagiellonian University,
ul.~Reymonta~4, PL-30-059 Krak\'ow, POLAND, 
e-mail: rams@im.uj.edu.pl

\end{document}